\documentclass[12pt,a4paper,reqno]{amsart}
\usepackage{amsmath,latexsym,epsfig,amsfonts,amssymb}
\usepackage{graphicx}

\newtheorem{theorem}{Theorem}
\newtheorem{lemma}[theorem]{Lemma}

\newcommand{\R}{\ensuremath{\mathbb{R}} }
\newcommand{\So}{\mathbb{S}}

\newcommand{\PLFip}{\ensuremath{\textup{PLF}^+([0,1])} }

\author{Nick Gill}
\author{Ian Short}
\email{Ian.Short@nuim.ie}
\date{\today}

\begin{document}

\begin{abstract}
Brin and Squier described how to classify the elements of the group of piecewise linear homeomorphisms of the real line according to conjugacy. We supply a short account of the key step in their classification. The results in this document are unoriginal. The proofs are a little sketchy in places.
\end{abstract}

\title{Piecewise linear homeomorphisms of the real line}
\maketitle

Let \PLFip denote the group of orientation preserving piecewise linear homeomorphisms of $[0,1]$ that are locally affine at all but a finite number of points in $[0,1]$. Let $\mathcal{F}$ denote the subset of $\PLFip$, which consists of homeomorphisms $f$ such that $f(x)>x$ for each $x\in(0,1)$. Membership of $\mathcal{F}$ is preserved by conjugation in $\PLFip$. Identifying the conjugacy classes that make up $\mathcal{F}$ is a key step in understand conjugacy in \PLFip and similar groups of piecewise linear homeomorphisms.

A \emph{node} of a map $f$ from \PLFip is a point $x$ in $(0,1)$ at which $f$ is not locally affine. Let $x_f$ denote the smallest node of $f$ in $(0,1)$. We also define $f'_-(x)$ and $f'_+(x)$ to be the left and right derivatives of $f$ at $x$. Define $f^*(x)=f'_+(x)/f'_-(x)$ so that $f^*(x)\neq 1$ if and only if $x$ is a node of $f$. By applying the chain-rule to the left and right derivatives of $f$, we deduce that there is a chain-rule for the operation $f\mapsto f^*$ which resembles the usual chain-rule for differentiation.

We associate to $f$ two objects $\alpha_f$ and $\beta_f$ that will 
be seen to be the only quantities necessary for a conjugacy classification in \PLFip of elements of $\mathcal{F}$. Define $\alpha_f=f'_-(x_f)$; this is the slope of the line segment in the graph of $f$ that protrudes from $0$. We can apply the chain-rule near $0$ to see that $\alpha$ is invariant under conjugation. To define $\beta_f$, first consider the the function $\phi_f:[x_f,f(x_f))\rightarrow \R^+$, given by
\begin{equation}\label{E: phi}
x\mapsto  \prod_{n=0}^\infty f^*(f^n(x)).
\end{equation}
Almost all terms $f^*(f^n(x))$ equal $1$. Hence, by the chain-rule, if $N$ is chosen to be sufficiently large then the image of $x$ in \eqref{E: phi} is equal to $(f^N)^*(x)$. The value of the product  remains unchanged if we extend to $n=-\infty$ because $f^*(x)=1$ for all $x\in(0,x_f)$. Let $\So$ denote the topological space $[0,1]$ with $0$ and $1$ identified. Let $\psi_f:\So\rightarrow\R^+$ be given by the equation
\begin{equation*}\label{E: psi}
\psi_f(s) = \phi_f(\alpha_f^sx_f).
\end{equation*}
 We consider the set $\mathcal{K}$ of functions from $\So$ to $\R^+$ that take the value $1$ at all but a finite number of points. We define an equivalence relation on $\mathcal{K}$ such that, for functions $u$ and $v$ in $\mathcal{K}$, $u$ is equivalent to $v$ if and only if there is a translation $t(x)=x+p \pmod 1$, $p\in\mathbb{R}$, such that $u=vt$. Finally, we define $\beta_f$ to be the equivalence class of $\psi_f$ in $\mathcal{K}$. That $\beta_f$ is invariant under conjugation in $\PLFip$ follows from the next two lemmas.

In both lemmas we use the notation $\phi_f^a$ to denote the map from $[a,f(a))$ to $\R^+$ given by \eqref{E: phi}. For $a\in (0,x_f)$, we define $\psi_f^a:\So\rightarrow\R^+$ by $\psi_f^a(s)=\phi_f^a(\alpha_f^s a)$.

\begin{lemma}\label{L: small}
If $a\in (0,x_f)$ then $\psi_f^a=\psi_f$.
\end{lemma}
\begin{proof}
Let $m=\alpha_f$. For a positive integer $n$, let $b=f^{-n}(x_f)$. For $s\in \So$,
\[
\psi^b_f(s)=\phi^b_f(m^s b)=\phi_f(f^n(m^sb))=\phi_f(m^sf^n(b))=\psi_f(s).
\]
If $a\in (b,f(b))$ we obtain
\[
\psi^a_f(s)=\psi_f(s+t),
\]
$t=\log_m (a/b)$, for $s\in\So$.
\end{proof}

\begin{lemma}\label{L: conj}
If $f$ and $g$ are two members of $\mathcal{F}$ that are conjugate in \PLFip then $\beta_f=\beta_g$.
\end{lemma}
\begin{proof}
Suppose that $h\in\PLFip$ and $hfh^{-1}=g$. We work in a small enough neighbourhood $U$ of $0$ that $f$, $g$, and $h$ are linear in $U\cup f(U)\cup g(U)$. Let $\lambda>0$ be such that $h(x)=\lambda x$ in $U$. Choose $a\in U$ and let $b=h(a)$. For $x\in[b,g(b))$ and suitably large integers $N$ we have,
 \[
(g^N)^*(x)=(hf^Nh^{-1})^*(x)=h^*(f^Nh^{-1}(x))(f^N)^*(h^{-1}(x)(h^{-1})^*(x) = (f^N)^*(x/\lambda).
\]
Let $m=\alpha_f$. Using Lemma~\ref{L: small} we obtain,
\[
\psi_g(s)=\phi_g^b(m^sb)=\phi_f^a(m^sb/\lambda)=\phi_f^a(m^sa)=\psi_f(s).
\]
\end{proof}

The remainder of this document is a proof of the following theorem.

\begin{theorem}\label{T: main}
Two maps $f,g\in\mathcal{F}$ are conjugate in \PLFip if and only if $\alpha_f=\alpha_g$ and $\beta_f=\beta_g$.
\end{theorem}

It remains to show that, for $f,g\in\mathcal{F}$, if $\alpha_f=\alpha_g$ and $\beta_f=\beta_g$ then $f$ and $g$ are conjugate in $\PLFip$. To prove this statement, we repeatedly conjugate by piecewise linear homeomorphisms with a single node. The result of this repeated conjugation will be a \emph{corner function}, that is, a member $f$ of $\mathcal{F}$ for which all nodes of $f$ occur in $[x_f,f(x_f))$. Given $p\in(0,1)$ and $0<\lambda<1$, there is a unique $h\in\PLFip$ with $h^*(p)=\lambda$ such that $h$ has only the single node $p$ in $(0,1)$. If $z_1,\dots,z_m$ are the nodes of $f$ in increasing order, then the nodes of $hfh^{-1}$ are contained within the set $\{h(z_1),\dots,h(z_m),h(p),hf^{-1}(p)\}$. Choose $p=z_m$ and $\lambda=f^*(z_m)$ then 
\[
(hfh^{-1})^*(h(p))=h^*(f(p))f^*(p)(h^{-1})^*(h(p))=f^*(p)(h^{-1})^*(h(p))=1.
\]
Thus $hfh^{-1}$ has nodes $h(z_1),\dots,h(z_{m-1}),hf^{-1}(z_m)$, and $(hfh^{-1})^*$ takes values \\$f^*(z_1),\dots,f^*(z_{m-1}),f^*(z_m)$ at these nodes. (Possibly $hf^{-1}(z_m)$ coincides with one of the $h(z_i)$.) We describe conjugation by $h$ as the \emph{elementary conjugation} of $f$. If $f$ is already a corner function and we apply the elementary conjugation, then the nodes of $f$ are permuted. (The last node of $f$ becomes the first node of $hfh^{-1}$ and all other nodes are shifted one place to the right. The values of the nodes remain unchanged.)

\begin{lemma}\label{L: corner}
Each member of $\mathcal{F}$ is conjugate to a corner function.
\end{lemma}
\begin{proof}
For a map $f\in\mathcal{F}$ that is not a corner function, let $x_f$  and $y_f$ be the smallest and largest nodes of $f$. Let $n\in\mathbb{N}$ be such that $y_f\in[f^n(x_f),f^{n+1}(x_f))$. Elementary conjugation of $f$ yields a function $g=hfh^{-1}$, and the node $y_f$ of $f$ is replaced by a node $hf^{-1}(y_f)$ of $g$ that lies in the interval $[g^{n-1}(x_g),g^n(x_g))$. Therefore repeatedly applying elementary conjugations yields a corner function. 
\end{proof}

\begin{lemma}\label{L: unique}
A corner function is uniquely specified by $\alpha$ and $\psi$.
\end{lemma}
\begin{proof}
Certainly a corner function $f$ is uniquely specified by $\alpha_f$ and $\phi_f$,  because $(x_f,\alpha_fx_f)$ is the first turning point in the graph of $f$, and all subsequent turning points are determined by $\phi_f$. If $f$ and $g$ are corner functions with $\alpha_f=\alpha_g$, $\psi_f=\psi_g$, but $x_f<x_g$, then $g$ must coincide with the function $x\mapsto \lambda f(\lambda^{-1}x)$, where $\lambda = x_g/x_f$. This function does not fix $1$, which is a contradiction. Therefore $x_f=x_g$. Hence $\phi_f=\phi_g$, hence $f=g$.\end{proof}

To complete the proof of Theorem~\ref{T: main}, suppose that $f$ and $g$ are maps in $\mathcal{F}$ such that $\alpha_f=\alpha_g$ and $\beta_f=\beta_g$. By Lemma~\ref{L: corner} we may assume that $f$ and $g$ are corner functions. We apply elementary conjugations to $f$, cycling the nodes of $f$, until $\psi_f=\psi_g$. From Lemma~\ref{L: unique} we deduce that $f=g$. This completes the proof of Theorem~\ref{T: main}.


\begin{thebibliography}{AW1}
\providecommand{\href}[2]{#2}
\bibitem{BrSq00} 
Brin, M. G. and Squier, C. C., 
Presentations, conjugacy, roots, and centralizers in groups of piecewise linear homeomorphisms of the real line,
 \emph{Comm. Algebra} (10) {\textbf 29} (2001), 4557--4596.
\end{thebibliography}
\end{document}